# THE NEW CLASSES OF THE GENETIC ALGORITHMS ARE DEFINED BY NONASSOCIATIVE GROUPOIDS.

S. SVERCHKOV, NOVOSIBIRSK STATE UNIVERSITY

**Abstract.** The genetic product of the groupoids, originating in the theory of DNA recombination, is introduced. It permits a natural generalization of the classical genetic algorithm. The full characterization of all three-element genetic groupoids gives an approach to construct the new classes of genetic algorithms. In the conclusion, we formulate some open problems in the theory of the genetic groupoids.

## 1. INTRODUCTION

By a *groupoid* $\mathbb{G} = (G, \cdot)$ we mean a non-empty set $G$ on which a binary (not necessary associative) operation " $\cdot$ " is defined. A finite set $A$ closed under the set of multiplications $M = \{\cdot_1, \cdot_2, ..., \cdot_n\}$ is called $n$-*groupoid* $\mathbb{A} = (A; \cdot_1, ..., \cdot_n)$ and the set $A$ is called the *solution space*.

Generally speaking any $n$-groupoid $\mathbb{A}$ defines $\mathbb{A}$-*algorithm*, which is an analog of classical genetic algorithm denoted throughout the paper by *GA*. Indeed, any operation of $M = \{\cdot_1, \cdot_2, ..., \cdot_n\}$ defines the *children*

$$a \cdot_i b, \ b \cdot_i a \in A, \qquad (1)$$

for any *parents* elements $a, b \in A$. Let us note that the operation (1) is a key process to obtain new generations in the iterative process of GA. To define the full procedure of an algorithm, we have to set up a fitness function for selection and the mutations. Let $f : A \to \mathbb{N}$ be an arbitrary function, which defines a fitness–based operator of selection the elements of any subset of $A$, and let $\varphi_i : A \to A$ be a finite set of the autoisomorphisms (mutations) of $A$. Then $\mathbb{A}$ – algorithm has the same scheme as GA with the replacement of crossovers by the operations of $M$ and with replacement the genetic mutation of GA by $\varphi_i$. Finally, to determine an $\mathbb{A}$-algorithm means simply the definition of multiplication table of $n$-*groupoid* $\mathbb{A} = (A; \cdot_1, ..., \cdot_n)$.

For example, in GA we have the solution space $A = \{(a_1, ..., a_{n+1}) \mid 0 \leq a_i \leq d_i\}$ and $n$ splicing operations (crossovers) $M = \{\times_1, \times_2, ..., \times_n\}$ on $A$ defined by rules

$$\forall a = (a_1, ..., a_{n+1}), b = (b_1, ..., b_{n+1}) \in A \quad (a_1, ..., a_{n+1}) \times_i (b_1, ..., b_{n+1}) = (a_1, ..., a_i, b_{i+1}, ..., b_{n+1}). \quad (2)$$

The splicing operations (2) define $n$-groupoid $\mathbb{GA}(n; d_1, d_2, ..., d_{n+1}) = (A; \times_1, ..., \times_n)$. We call $\mathbb{GA}(n; d_1, d_2, ..., d_{n+1})$ the *splicing* $n$-groupoid. Therefore classical genetic algorithm is a $\mathbb{GA}(n; d_1, d_2, ..., d_{n+1})$-algorithm. It is easy to check that the splicing operations are associative (see [1] for more details of splicing operations), i.e.

$$\forall a, b, c \in A, \ 1 \leq i \leq n-1 \quad (a \times_i b) \times_i c = a \times_i (b \times_i c).$$

However, a "mixture" of splicing operations is not associative. Indeed, we have

$$((111) \times_2 (101)) \times_1 (000) = (111) \times_1 (000) = (100),$$
$$(111) \times_2 ((101) \times_1 (000)) = (111) \times_2 (100) = (110),$$

and

$$((111) \times_2 (101)) \times_1 (000) \neq (111) \times_2 ((101) \times_1 (000)).$$

We have to indicate two more important properties of the splicing $n$-groupoids $\mathbb{GA}(n; d_1, d_2, ..., d_{n+1})$. For every splicing operation $\cdot \in M$, we have

$$\forall a \in A \quad a \cdot a = a, \tag{3}$$
$$\forall a, b \in A \quad a \cdot b = b \cdot a \implies a = b. \tag{4}$$

Those natural restrictions of the operations eliminate the repetitions in the genetic algorithm process.

We call a groupoid ($n$-groupoid) *idempotent* if it satisfies (3); *nowhere commutative* (see [2] for the definition) if it satisfies (4); and *genetic* if it satisfies (3) and (4). We will denote by $\mathfrak{Gen}(n)$ the class of all genetic $n$-groupoids, for $n \in \mathbb{N}$; by $\overline{\mathfrak{Gen}} = \bigcup_{n \in \mathbb{N}} \mathfrak{Gen}(n)$ the class of all genetic $n$-groupoids, for all $n \in \mathbb{N}$. Set $\mathfrak{Gen} = \mathfrak{Gen}(1)$. Consequently, every groupoid $(A, \times_i) \in \mathfrak{Gen}$ is a *genetic* semigroup. One can easily check that $\mathbb{GA}(n; d_1, d_2, ..., d_{n+1})$ is a genetic $n$-groupoid.

Naturally, we have to consider only the genetic $n$-groupoids on a way of the construction of new classes of genetic algorithms. Even so, an $\mathbb{A}$ – algorithm for a genetic $n$-groupoid of a large enough solution space has useless meaning for the computer applications. The multiplication table of such $n$-groupoids is too large to operate with. The obvious advantage of GA is simplest and easy computer applicable multiplication table (2) of splicing $n$-groupoids. Therefore, our goal is to invent a construction of generating the subclasses of $\overline{\mathfrak{Gen}}$ which preserves those benefits of splicing $n$-groupoids. It will give us an approach to construct the new classes of genetic algorithms, which are computable ones.

2. GENETIC PRODUCT

Let $\mathbb{A} = (A, M)$, $\mathbb{B} = (B, N)$ be a $n$-groupoid with the solution space $A = \{a_i\}$, and with the operations $M = \{\cdot_1, \cdot_2, ...., \cdot_n\}$; and a $m$-groupoid with the solution space $B = \{b_j\}$, and with operations $N = \{*_1, *_2, ..., *_m\}$, respectively. Let us define the $(n+m+1)$ operations on the set $A \times B = \{(a, b) \mid a \in A, b \in B\}$ by the next rules:

$$(a,b) \cdot_x (c,d) = \begin{cases} (a \cdot_x c, d), & x \leq n, \\ (a, d), & x = n+1, \\ (a, b *_{x-n-1} d), & n+2 < x \leq n+m+1. \end{cases} \tag{5}$$

for any $(a,b),(c,d) \in A \times B$. The obtained $(n+m+1)$-groupoid $\mathbb{A} *_G \mathbb{B}$ with the operation set $\overline{M} = \{\cdot_1, \cdot_2, ..., \cdot_{n+m+1}\}$ is called the *genetic product* of $\mathbb{A}$ and $\mathbb{B}$.

It follows immediately that genetic product preserves genetic property of the groupoids, that is

$$\mathbb{A} \in \mathfrak{Gen}(n), \mathbb{B} \in \mathfrak{Gen}(m) \Rightarrow \mathbb{A} \times_G \mathbb{B} \in \mathfrak{Gen}(n+m+1).$$

Set $\mathcal{GA}(n,d) = \mathcal{GA}(n; \underbrace{d,...,d}_{n+1})$. The trivial case $G(0,d)$ means simply the set $\{0,1,...,d\}$ without any operations. The genetic product of $\mathbb{A} = (A, M)$ and $G(0,d)$ is called the *genetic extension* of $\mathbb{A}$ and it is denoted by $G(\mathbb{A}) = \mathbb{A} *_G G(0,d)$. It is clear from (5) that solution space $G(A) = \{(a, j) \mid a \in A, \ 0 \le j \le d\}$ defines $(n+1)$-groupoid $G(\mathbb{A})$ with the multiplication table:

$$(a, j) \cdot_x (b, t) = \begin{cases} (a \cdot_x b, t), & 1 \le x \le n, \\ (a, t), & x = (n+1). \end{cases}$$

The straightforward proof provides us with useful properties of the genetic product.

**Lemma 1.**
$$\forall \mathbb{A}, \mathbb{B}, \mathbb{C} \in \overline{\mathfrak{Gen}} \quad (\mathbb{A} *_G \mathbb{B}) *_G \mathbb{C} \cong \mathbb{A} *_G (\mathbb{B} *_G \mathbb{C}), \tag{6}$$

and for the general case, we have $\mathbb{A} *_G \mathbb{B} \not\cong \mathbb{B} *_G \mathbb{A}$.

$\mathfrak{G} = (\mathbb{A} \in \mathfrak{Gen}, *_G)$ is a semigroup, and besides $\forall \mathbb{A} \in \mathfrak{G} \ \mathbb{A} \in \overline{\mathfrak{Gen}}$.

$$\forall n, m \in \mathbb{N} \quad GA(n+m+1; a,b) \cong GA(n,a) *_G GA(m,b),$$

$$GA(1,d) \cong G(G(0,d)) \cong G(0,d) \times_G G(0,d), \ GA(n,d) \cong \underbrace{G(...(G(1,d)...)}_{n},$$

$$G(2n,d) \cong \underbrace{G(1,d) \times_G ... \times_G G(1,d)}_{n-1}, \quad G(2n+1,d) \cong G(2n,d) \times_G G(0,d).$$

Let $\mathbb{A} = \mathbb{A}_1 \times_G ... \times_G \mathbb{A}_n$, $\mathbb{A}_i \in \overline{\mathfrak{Gen}}$, and $\varphi : \mathbb{A}_i \to \mathbb{A}_i$ be an autoisomorphism of $\mathbb{A}_i$, then $\varphi_G(a_1,...,a_n) = (a_1,...,\varphi(a_i),...,a_n)$ is an autoisomorphism of $\varphi_G : \mathbb{A} \to \mathbb{A}$.

Conversely, suppose $\phi(a_1,...,a_n) = (\phi(a_1),...,\phi(a_n))$ is an autoisomorphism of $\mathbb{A}$, then $\phi_i(a_i) = \phi(a_i)$ is an autoisomorphism of $\mathbb{A}_i$, $1 \le i \le n$.

The autoisomorphism group $Aut(\mathbb{A})$ of $\mathbb{A}$ is fully defined by the autoisomorphism groups $Aut(\mathbb{A}_i)$ of $\mathbb{A}_i$, $1 \le i \le n$

We call $\mathbb{A} = \mathbb{A}_1 \times_G ... \times_G \mathbb{A}_n$, $\mathbb{A}_i \in \overline{\mathfrak{Gen}}$ the *genetic decomposition* of $\mathbb{A}$ into $\mathbb{A}_i \in \overline{\mathfrak{Gen}}, 1 \leq i \leq n$.

## 3. MAIN RESULTS

Let us note that the genetic product $\mathbb{A} = \mathbb{A}_1 \times_G ... \times_G \mathbb{A}_n = (\mathbb{A}_1 \times ... \times \mathbb{A}_n; *_i)$ of the genetic groupoids $\mathbb{A}_i = (A_i, \cdot_i)$, $1 \leq i \leq n$ provides the simple computer applicable multiplication rules by modulo of the multiplications of $\mathbb{A}_i = (A_i, \cdot_i)$. Indeed, from (5), (6) we have

$$\forall a_i, b_i \in A_i \quad (a_1,...,a_n) *_i (b_1,...,b_n) = \begin{cases} (a_1,...,a_{i-1}, b_i,...,b_n), & i = 2k, \\ (a_1,...,(a_i \cdot_i b_i),...,b_n), & i = 2k+1. \end{cases}$$

From this we conclude that the complexity of computer calculations of $\mathbb{A} \cong \mathbb{A}_1 \times_G ... \times_G \mathbb{A}_n$-algorithm for small dimensions of the groupoids $\mathbb{A}_i$ (let us say $|\mathbb{A}_i| \leq 4$) is no substantially greater than the complexity of classical GA (respectively, complexity of $GA(2n, 2)$). Moreover, the genetic decomposition of $\mathbb{A} = \mathbb{A}_1 \times_G ... \times_G \mathbb{A}_n$ transfers us from $\overline{\mathfrak{Gen}}$ into $\mathfrak{Gen}$ to operate with single operation groupoids. And moreover, by lemma 1, all mutations of $\mathbb{A} = \mathbb{A}_1 \times_G ... \times_G \mathbb{A}_n$ is defined by the autoisomorphisms of $\mathbb{A}_i$.

Let $\mathbb{A} = \mathbb{A}_1 \times_G ... \times_G \mathbb{A}_n$ be the genetic product of arbitrary finite genetic semigroups (associative genetic groupoids) $\mathbb{A}_i = (A_i, \cdot_i)$. We first prove that $\mathbb{A}$-algorithm is a classical genetic algorithm for a splicing $(2n-1)$-groupoid.

The semigroup $GA(1; n, m)$ is called the *rectangular band* of $\{1,...,n\} \times \{1,...,m\}$. Obviously, this semigroup has a multiplication table $(x,y) \cdot (a,b) = (x,b)$. For more details of the rectangular bands we refer the readers to [2, 3]. Let us denote by $\mathfrak{RBand}$ the class of the rectangular bands. It is well known that any genetic semigroup is a rectangular band, and the class $\mathfrak{RBand}$ is a variety defined by two identities

$$\begin{cases} \forall a,b,c & (ab)c = a(bc), \\ \forall a,b & (ab)a = a. \end{cases} \quad (7)$$

On account of the above remark, we obtain a characterization of a genetic decomposition of $\mathbb{A}$ into the associative genetic groupoids (genetic semigroups).

**Theorem 1.** Let $\mathbb{A}_i = (A_i, \cdot_i)$, $1 \leq i \leq n$ be a family of finite genetic semigroups then $\mathbb{A} = \mathbb{A}_1 \times_G ... \times_G \mathbb{A}_n \cong G(2n-1; d_1,...,d_{2n})$ for some $d_i \in \mathbb{N}$.

Consequently, there are no new classes of genetic algorithms in associative case. Therefore, to obtain a generalization of GA if it any exists, we have to describe nonassociative genetic groupoids of the small dimensions. The procedure (1) of creating new generations in $\mathbb{A}$-algorithm shows that it is sufficient to characterize groupoids up to isomorphism and anti-isomorphism. We will do it for three–element ones.

It is easy to check that any two element genetic groupoid is a semigroup, i.e. they are in the variety $\mathfrak{RBand}$. It follows that we have to deal with at least 3-element groupoids. By a computer study, all non-isomorphic 3-element groupoids were described in [4]. There are 3330 of them. It is too much even to observe them. Fortunately, the class $\mathfrak{Gen}_3$ of genetic (not isomorphic or anti- isomorphic) 3-element groupoids consists only of 18 groupoids. A characterization of $\mathfrak{Gen}_3$ can be easily made by "hand" calculations.

Let us denote by $\begin{pmatrix} i & j & k \\ x & y & z \end{pmatrix}$, $i \neq x, j \neq y, k \neq z$ the genetic groupoid $A = (\{a_1, a_2, a_3\}, \cdot)$ with multiplication table $\begin{pmatrix} a_1 & a_i & a_j \\ a_x & a_2 & a_k \\ a_y & a_z & a_3 \end{pmatrix}$.

Directly applying only three automorphisms (permutations) $\varphi_{ij} : a_i \to a_j$ to the multiplication table of $A = (\{a_1, a_2, a_3\}, \cdot)$ we obtain the description of $\mathfrak{Gen}_3$.

Set $\mathfrak{M} = \{A_i, 1 \leq i \leq 17\} =$

$$= \left\{ \begin{pmatrix} 000 \\ 111 \end{pmatrix}, \begin{pmatrix} 000 \\ 222 \end{pmatrix}, \begin{pmatrix} 000 \\ 112 \end{pmatrix}, \begin{pmatrix} 000 \\ 121 \end{pmatrix}, \begin{pmatrix} 000 \\ 211 \end{pmatrix}, \begin{pmatrix} 000 \\ 221 \end{pmatrix}, \begin{pmatrix} 000 \\ 212 \end{pmatrix}, \begin{pmatrix} 000 \\ 122 \end{pmatrix}, \begin{pmatrix} 100 \\ 221 \end{pmatrix}, \right.$$

$$\left. \begin{pmatrix} 011 \\ 122 \end{pmatrix}, \begin{pmatrix} 011 \\ 221 \end{pmatrix}, \begin{pmatrix} 002 \\ 121 \end{pmatrix}, \begin{pmatrix} 020 \\ 112 \end{pmatrix}, \begin{pmatrix} 200 \\ 112 \end{pmatrix}, \begin{pmatrix} 111 \\ 020 \end{pmatrix}, \begin{pmatrix} 111 \\ 200 \end{pmatrix}, \begin{pmatrix} 111 \\ 002 \end{pmatrix} \right\}.$$

**Theorem 2.** $\mathfrak{Gen}_3 = \mathfrak{M} \cup \left\{ \begin{pmatrix} 001 \\ 122 \end{pmatrix} \right\}$.

The class $\mathfrak{M}$ consists of all three–element genetic nonassociative groupoids up to isomorphism or anti- isomorphism. There is only one three–element genetic associative groupoid in $\mathfrak{Gen}_3$. It is a semigroup $\begin{pmatrix} 001 \\ 122 \end{pmatrix} \in \mathfrak{RBand}$, and it is precisely the rectangular band of $\{1\} \times \{1, 2, 3\}$.

Application of Theorem 2 gives us the new classes of computer applicable $\mathbb{A} = \mathbb{A}_1 \times_G ... \times_G \mathbb{A}_m$-algorithms. They are strictly nonassociative $\mathbb{A}$-algorithms for $A_i \in \mathfrak{M}$, $m \in \mathbb{N}$, and partially associative $\mathbb{A}$ - algorithms for $A_i \in \mathfrak{M} \bigcup_{d \in \mathbb{N}} \{G(0, d)\}$.

### 4. CONCLUSION

It is clear that the class $\mathfrak{Gen}$ of all genetic groupoids defines all natural generalizations of the classical genetic algorithms, by modulo of the genetic decomposition. Therefore, any characterization of $\mathfrak{Gen}$ may be an extremely useful. One can check that $\mathfrak{Gen}$ forms a quasi-variety defined by the identity (3) and quasi-identity (4).

© 1990 Sverchkov S.R.

It is not hard to see that $\mathfrak{Gen}$ does not form a variety. It is closed under formation of subgroupoids and direct products, but not of homomorphic images. The homomorphic image of genetic groupoid $\begin{pmatrix} 000 \\ 111 \end{pmatrix}$, defined by $a_0 = a_1$, is not genetic. From this it follows that the class $\mathfrak{Gen}$ can't be defined by any set of identities. But there are some subvarieties in $\mathfrak{Gen}$. For example $\mathfrak{RBand} \subset \mathfrak{Gen}$, by (7).

From the above remarks, one may set up some problems, whose answers aren't known to author.

Problem 1. Describe maximal subvarieties of $\mathfrak{Gen}$.

The variety $\mathfrak{RBand}$ has an elegant and complete description of the lattice of the all subvarieties (see [5] for more details). In particular, it is countable, complete, and distributive.

Problem 2. Describe the lattice of the all subvarieties of $\mathfrak{Gen}$.

S. SVERCHKOV,
NOVOSIBIRSK STATE UNIVERSITY,
e-mail: sverchkovsr@yandex.ru